\documentclass[number,citesort,seceqn,dvips]{arxbj}
\usepackage{upgreek}

\aid{0}
\volume{18}
\issue{1}
\pubyear{2012}
\firstpage{119}
\lastpage{136}
\doi{10.3150/10-BEJ325}

\makeatletter
\newcommand{\reel}{\mathbb{R}}
\newtheorem{propo}{Proposition}[section]
\newtheorem{theo}[propo]{Theorem}
\newremark{rem}[propo]{Remark}
\newtheorem{lemma}[propo]{Lemma}
\newproclaim{defini}[propo]{Definition}
\makeatother

\begin{document}
\begin{frontmatter}

\title{Independence properties of the Matsumoto--Yor type}
\runtitle{Independence properties of the Matsumoto--Yor type}

\begin{aug}
\author{\fnms{A.E.} \snm{Koudou}\corref{}\thanksref{e1}\ead[label=e1,mark]{Efoevi.Koudou@iecn.u-nancy.fr}} \and
\author{\fnms{P.} \snm{Vallois}\thanksref{e2}\ead[label=e2,mark]{Pierre.Vallois@iecn.u-nancy.fr}}
\runauthor{A.E. Koudou and P. Vallois}
\address{Institut Elie Cartan, Laboratoire de Math\'{e}matiques, B.P.
239, F-54506 Vandoeuvre-l\`{e}s-Nancy CEDEX, France. \printead{e1,e2}}
\end{aug}

\received{\smonth{5} \syear{2009}}
\revised{\smonth{8} \syear{2010}}

%
\begin{abstract}
We define Letac--Wesolowski--Matsumoto--Yor (LWMY) functions as
decreasing functions from $(0,\infty)$ onto $(0,\infty)$ with the
following property: there exist independent, positive random variables
$X$ and $Y$ such that the variables $f(X+Y)$ and $f(X)-f(X+Y)$ are
independent. We prove that, under additional assumptions, there are
essentially four such functions. The first one is $f(x)=1/x$. In this
case, referred to in the literature as the \textit{Matsumoto--Yor
property}, the law of $X$ is generalized inverse Gaussian while $Y$ is
gamma distributed. In the three other cases, the associated densities
are provided. As a consequence, we obtain a new relation of convolution
involving gamma distributions and Kummer distributions of
type~2.\looseness=1
\end{abstract}

%

\vspace*{12pt}
\end{frontmatter}
%

\section{Introduction}

Many papers have been devoted to generalized inverse Gaussian
($\operatorname{GIG}$) distributions since their definition by Good
\cite{good} (see, e.g., \cite{barnd,letacsesh,val1,val2}).

 The $\operatorname{GIG}$ distribution with parameters $\mu\in\reel$,
$a,b>0 $ is the probability measure
%
\begin{equation} \label{defgig}
\operatorname{GIG}(\mu, a, b)(\mathrm{d}x)= \biggl( \frac{b}{a}
\biggr)^\mu\frac{x^{\mu-1}}{ 2K_\mu(ab)}
\mathrm{e}^{-(a^2x^{-1}+b^2x)/2}{\mathbf 1}_{(0,\infty)}(x)\,\mathrm{d}x,
\end{equation}
where $K_\mu$ is the classical McDonald special
function.



\begin{longlist}
\item[\textbf{(1)}] We stress the close links between $\operatorname{GIG}$, gamma
distributions and the function $f_0(x)=1/x$ ($x>0$).

\begin{longlist}
\item[\textbf{(a)}] The family of $\operatorname{GIG}$ distributions is invariant
under $f_0$: we can easily deduce from~(\ref{defgig}) that the image of
$\operatorname{GIG}(\mu, a, b)$ by $f_0$ is $\operatorname{GIG}(-\mu,
b, a)$.

\item[\textbf{(b)}] Barndorff-Nielsen and Halgreen \cite{barnd} proved that
%
\begin{equation} \label{convogig}
\operatorname{GIG}(-\mu,a,b) * \gamma\biggl(\mu, \frac{b^2}{2}\biggr)
=\operatorname{GIG}(\mu,a,b),\qquad   \mu, a,
 b>0,
\end{equation}
where $\gamma(\mu, b^2 /2)(\mathrm{d}x)= \frac{b^{2\mu}}{2^{\mu}\Gamma(\mu
)} x^{\mu-1}
\exp-\frac{b^2}{2}x\mathbf{1}_{(0,\infty)}(x)\,\mathrm{d}x$.\vadjust{\goodbreak}

Therefore, if $ X \sim\operatorname{GIG}(-\lambda, a, a)$ and
$Y\sim\gamma (\lambda, a^2 /2)$ are independent random variables, then
%
\begin{equation} \label{equiv}
X \stackrel{(d)}{=}f_0(X+Y).
\end{equation}
Letac and Seshadri \cite{letacsesh} proved that (\ref{equiv}) characterizes
$\operatorname{GIG}$ distributions of the type
$\operatorname{GIG}(-\lambda, a, a)$.

\item[\textbf{(c)}]
Almost sure realizations of (\ref{convogig}) have been given by
Bhattacharya and Waymire~\cite{bhat} in the case $\mu=\frac{1}{2}$,
Vallois \cite{val2} for any $\mu>0$ by means of a family of transient
diffusions and
Vallois~\cite{val1}, theorem on page 446, in terms of random walks.
\end{longlist}
\item[\textbf{(2)}]
The so-called \textit{Matsumoto--Yor property} is the following:
let $X$ and $Y$ be two independent random variables such that
%
\begin{equation} \label{lawsxy}
X \sim\operatorname{GIG}(-\mu, a, b),\qquad   Y \sim\gamma(\mu, b^2 /2),\qquad
(\mu, a, b>0).
\end{equation}
Then,
%
\begin{equation} \label{defininguv}
U:= \frac{1}{X+Y}=f_0(X+Y),\qquad   V:= \frac{1}{X}-
\frac{1}{X+Y}=f_0(X)-f_0(X+Y)
\end{equation}
are independent and
%
\begin{equation} \label{lawsuv}
U \sim\operatorname{GIG}(-\mu, b, a),\qquad   V \sim\gamma(\mu, a^2 /2).
\end{equation}
\end{longlist}

The case $a=b$ was proven by Matsumoto and Yor \cite{mats} and a nice
interpretation of this property via Brownian motion was given by
Matsumoto and Yor \cite{mats2}. The case $\mu= -\frac{1}{2}$ of the
Matsumoto--Yor property can be retrieved from an independence property
established by Barndorff-Nielsen and Koudou \cite{barn} (see \cite{koud}).


Letac and Wesolowski \cite{letacwes} proved that the Matsumoto--Yor
property holds for any~$\mu,\allowbreak a, b >0$ and characterizes the
$\operatorname{GIG}$ distributions. More precisely, consider two
independent and non-Dirac positive random variables $X$ and $Y$ such
that $U$ and $V$ defined by~(\ref{defininguv}) are independent. There
then exist $\mu, a,b>0$ such that (\ref{lawsxy}) holds.

The starting point of this paper is to study the link between the function
$f_0\dvtx x \mapsto1/x$ and the $\operatorname{GIG}$ distributions in the
Matsumoto--Yor
property.

Obviously, the Matsumoto--Yor property can be re-expressed as follows:
the image of the probability measure (on $\reel_+ ^2$)
$\operatorname{GIG}(-\mu, a,b) \otimes\gamma(\mu,b^2/2)$ by the
transformation $T_{f_0}\dvtx (x,y) \mapsto (f_0(x+y), f_0(x)-f_0(x+y))$ is
the probability measure $\operatorname{GIG}(-\mu, b,a)
\otimes\gamma(\mu,a^2/2)$. This formulation of the Matsumoto--Yor
property, joined with the Letac and Wesolowski result, leads us to
determine the triplets $(\mu_X, \mu_Y, f)$ such that:

\begin{longlist}
\item[(a)] $\mu_X, \mu_Y$ are probability measures on
$(0,\infty)$;\vadjust{\goodbreak}

\item[(b)] $f\dvtx  (0,\infty) \rightarrow(0,\infty)$ is bijective and decreasing;

\item[(c)] if $X$ and $Y$ are independent random variables such that $X\sim
\mu_X$ and $Y\sim\mu_Y$, then the random variables $U=f(X+Y)$ and
$V=f(X)-f(X+Y)$ are independent.
\end{longlist}

Unfortunately, we have not been able to solve this question without
restriction. Our method can be applied provided that $f$ is smooth and
$\mu_X$ and $\mu_Y$ have smooth density functions (see Theorem
\ref{cns} for details). After long and sometimes tedious calculations,
we prove (see Theorem~\ref{theoclass}) that there are only four
classes, ${\mathcal F}_1,\ldots,{\mathcal F}_4$, of functions~$f$ such that
$T_f$ keeps the independence property. Then, for any $f \in{\mathcal F}_i ,
1\leq i \leq4,$ we have been able to give the corresponding
distributions of $X$ and $Y$ and the related laws of $U$ and $V$ (for
${\mathcal F}_2, {\mathcal F}_3$ and ${\mathcal F}_4,$ see Theorems~\ref{casfun},
\ref{casfdelta} and Remark~\ref{newrem}). The first class, ${\mathcal F}_1
= \{ \alpha/x;  \alpha>0 \},$ corresponds to the known case $f=f_0$.
This case, as mentioned in Remark~\ref{casyor}, allows us recover,
under stronger assumptions, the result of Letac and Wesolowski that the
only possible distributions for $X$ and $Y$ are $\operatorname{GIG}$
and gamma, respectively. The proof of Letac and Wesolowski is
completely different from ours since the authors make use of Laplace
transforms and a characterization of the $\operatorname{GIG}$ laws as
the distribution of a continued fraction with gamma entries.
We have not been able to develop a proof as elegant as theirs because,
with $f=f_0,$ we have algebraic properties (e.g., continued fractions),
while these properties are lost if we start with a general
function~$f$.\looseness=1

It is worth pointing out that one interesting feature of our analysis
is an original characterization of the families of distributions $\{
\beta_\alpha(a,b,c);  a, b, \alpha> 0,  c\in\reel\}$ and the Kummer
distributions $\{ K^{(2)}(a,b,c);  a, c > 0, b \in\reel\} $ (see
(\ref{kummer}) and (\ref{betaalpha}), respectively). The Kummer
distributions appear as the laws of some random continued fractions (see~\cite{mark}, page~3393, mentioning a~work by Dyson \cite{Dys} in
the setting of random matrices).

As by-products of our study, we obtain new relations for convolution.
For simplicity, we only detail the case of Kummer distributions of type
2:
%
\begin{equation} \label{convokummer}
K^{(2)}(a,b,c) * \gamma(b,c)=
K^{(2)}(a+b,-b,c).
\end{equation}
Obviously, this relation is similar to (\ref{convogig}).

Inspired by the result of Letac and Wesolowski \cite{letacwes} and Theorem
\ref{cas1}, we can ask (for the purposes of future research) whether a
characterization of Kummer distributions could be obtained via an
``algebraic'' method.

As recalled in the above item (c), there are various almost sure
realizations of (\ref{convogig}) and of the convolution coming from the
Matsumoto--Yor property. One interesting open question derived from our
study would be to determine a random variable $Z$ with distribution $
K^{(2)}(a+b,-b,c)$ which can be decomposed as the sum of two explicit
independent random variables $X$ and $Y$ such that $X\sim
K^{(2)}(a+b,-b,c)$ and $Y\sim\gamma(b,c)$.

The paper is organized as follows. We state our main results in
Section \ref{sec2}. In Section \ref{sec3} we give a key differential equation involving
$f$ and the log densities of the independent random variables $X$ and
$Y$ such that $f(X+Y)$ and $f(X)-f(X+Y)$ are independent (see Theorem~\ref{cns}).
Based on this equation, we prove (see Theorem~\ref{classificationFbis}) that there are only four classes of such
functions~$f$. The theorems stated in Section \ref{sec2} are proved in Section
\ref{sec4}; however, one technical proof has been postponed to the \hyperref[appendix]{Appendix}.

\section{Main results}\label{sec2}
\setcounter{equation}{7}
\begin{defini}
Let $f\dvtx (0,\infty) \rightarrow(0,\infty)$ be a decreasing and
bijective function.
\begin{longlist}
\item[(1)] We consider the transformation associated with $f$
\begin{eqnarray}\label{transformation}
T_f\dvtx  (0,\infty)^2 &\rightarrow& (0,\infty)^2,\nonumber\\[2pt]
(x,y) &\mapsto&
\bigl(f(x+y),f(x)-f(x+y)\bigr).
\end{eqnarray}
The transformation $T_f$ is one-to-one and if $f^{-1}$ is the inverse
of $f$, then
%
\begin{equation} \label{transfinverse}
(T_f)^{-1}= T_{f^{-1}}.
\end{equation}

\item[(2)] Let $X$ and $Y$ be two independent and positive random variables.
Let us define
%
\begin{equation} \label{defuv}
(U,V)= T_f(X,Y) = \bigl(f(X+Y),f(X)-f(X+Y)\bigr).
\end{equation}
$f$ is said to
be an \textup{LWMY function with respect to} $(X,Y)$ if the random
variables $U$ and~$V$ are independent.
$f$ is said to be an LWMY function if it is an LWMY function with
respect to some random vector $(X,Y)$.
\end{longlist}
\end{defini}

One aim of this paper is to characterize LWMY functions.
Let us introduce
%
\begin{eqnarray} \label{expo}
f_1(x)&=&\frac{1}{e^x-1},\qquad x>0,\\[2pt] \label{expoinv}
g_1(x)&=&f_1^{-1}(x)=\ln\biggl( \frac{1+x}{x} \biggr),\qquad   x>0
\end{eqnarray}
and, for $\delta>0$,
%
\begin{equation} \label{involution}
f_\delta^{*}(x)=\log\biggl( \frac{e^x+\delta-1}{e^x-1} \biggr),\qquad   x>0.
\end{equation}

\begin{theo} \label{theoclass}
Let $f\dvtx(0,\infty) \rightarrow(0,\infty)$ be decreasing and
bijective. Under some additional assumptions (see Theorem
\ref{theocns}, (\ref{cond2}) and (\ref{serie})), $f$ is an LWMY
function if and only if either $f(x)=\frac{\alpha}{x}$, $
f(x)=\frac{1}{\alpha}f_1(\beta x)$, $ f(x)=\frac{1}{\alpha
}g_1(\beta
x)$ or $ f(x)=\frac{1}{\alpha}f_\delta^* (\beta x)$ for some $\alpha,
\beta$, $\delta>0$.
\end{theo}

\begin{rem} \label{remark1}
(1) The four classes of LWMY functions are ${\mathcal F}_1= \{ \alpha/x;
\alpha>0\}, \ {\mathcal F}_2= \{ \frac{1}{\alpha} f_1(\beta x);
\alpha,\beta>0\},$ $ {\mathcal F}_3= \{ \frac{1}{\alpha} g_1(\beta x);
\alpha,\beta>0\}$ and ${\mathcal F}_4= \{ \frac{1}{\alpha} f_\delta^*(\beta
x); \alpha
,\beta>0\}.$\vspace*{-6pt}
\begin{longlist}
\item[(2)] It is clear that if $f$ is an LWMY function, then the functions
$f^{-1}$ and $x\mapsto\frac{1}{\alpha}f(\beta x), \alpha, \beta
>0,$ are LWMY functions. 

\item[(3)] The image of ${\mathcal F}_2$ by the map $f \mapsto f^{-1}$ is ${\mathcal
F}_3$. The functions $x\mapsto\alpha/ x$ and $f_\delta$ are
involutive.
\end{longlist}

In the sequel, we focus on the three new cases: either $f=f_1$, $f=g_1$
or $f=f_\delta^*$ and in each case, we determine the laws of the
related random variables.
\end{rem}

\subsection{The cases $f=g_1$ and $f=f_1$}\vspace*{-3pt}

\hspace*{9pt}\textbf{(a)} Recall the definitions of the gamma distribution $\gamma
(\lambda, c)(\mathrm{d}x)=\break \frac{c^\lambda}{\Gamma(\lambda)} x^{\lambda-1}
\mathrm{e}^{-cx} {\mathbf1}_{(0,\infty)}(x)\,\mathrm{d}x    (\lambda, c>0)$ and the beta
distribution $\mathrm{Beta} (a,b) (\mathrm{d}x)=\break\frac{\Gamma(a+b)}{\Gamma(a)
\Gamma(b)} x^{a-1} (1-x)^{b-1} {\mathbf1}_{ \{0<x<1\} }\,\mathrm{d}x\,
   (a, b>0).$
Consider (see, e.g., \cite{kotz}, or \cite{Nagupt} and
the references therein) the \textit{Kummer distribution of type 2}:\vspace*{-2pt}
%
\begin{equation} \label{kummer}
K^{(2)}(a,b,c): =\alpha(a, b,c) x^{a -1} (1+x)^{-a-b} \mathrm{e}^{-cx} {\mathbf
1}_{(0,\infty)}(x)\,\mathrm{d}x,\qquad  a,c>0,  b \in\reel,\vspace*{-2pt}
\end{equation}
where $\alpha(a, b, c)$ is a normalizing constant.

Associated with a couple $(X,Y)$ of positive random variables,
consider\vspace*{-2pt}
%
\begin{equation} \label{uvenfxy}
(U,V):= T_{f_1}(X,Y) = \biggl( \frac{1}{\mathrm{e}^{X+Y}-1}, \frac{1}{\mathrm{e}^{X}-1} -
\frac{1}{\mathrm{e}^{X+Y}-1} \biggr).\vspace*{-2pt}
\end{equation}

In Theorems \ref{casfun} and \ref{cas1} below, we suppose that all
random variables have positive and twice differentiable densities.

First, we consider the case $f=f_1$. We determine the distributions of
$X$ and $Y$ such that $f_1$ is an LWMY function associated with
$(X,Y)$.\vspace*{-3pt}
\begin{theo} \label{casfun}
\textup{(1)} Consider two positive and independent random variables $X$ and $Y$.
The random variables $U$ and $V$ defined by (\ref{uvenfxy}) are
independent if and only if
the densities of $Y$ and $X$ are, respectively,\vspace*{-2pt}
%
\begin{eqnarray} \label{sol35}
p_Y(y)& = & \frac{\Gamma(a+b)}{\Gamma(a) \Gamma(b)} (1-\mathrm{e}^{-y})^{b-1}
\mathrm{e}^{-a y} {\mathbf1}_{ \{y>0\} }, \\[-2pt] \label{sol36}
p_X(x) & = & \alpha(a+b, c, -a) \mathrm{e}^{-(a+b ) x} (1-\mathrm{e}^{-x})^{-b
-1}\nonumber\\ [-9pt]\\ [-9pt]
& &{}\times \exp\biggl( -c \frac{\mathrm{e}^{-x}}{1-\mathrm{e}^{-x}} \biggr) {\mathbf1}_{ \{x>0\} },\nonumber\vspace*{-2pt}
\end{eqnarray}
where $a$, $b$ and $c$ are constants such that $a,b,c>0$ and $\alpha
(a+b, c, -a)$ is the constant from equation (\ref{kummer}). Thus, the
law of $Y$ is the image of the $\operatorname{Beta}(a,b)$ distribution by the
transformation $z\in(0,1)\mapsto-\log z \in(0,\infty)$, while the
law of the variable $f_1(X)$ is $K^{(2)}(a+b,-b,c)$ (see equation
(\ref{kummer})).

\textup{(2)} If \textup{(1)} holds, then $U \sim K^{(2)}(a,b,c)$ and $ V \sim\gamma(b,c).$\vspace*{-3pt}
\end{theo}

The proof of Theorem \ref{casfun} will be given in Section \ref{sec4}.\vspace*{-3pt}
\begin{rem} \label{newrem}
Since $g_1=f_1^{-1}$, Remark \ref{remark1} and Theorem \ref{casfun}
imply that the random variables associated with the LWMY function
$g_1$ are the random variables $U$ and $V$ distributed as in item~2 of
Theorem \ref{casfun}.\vspace*{-3pt}
\end{rem}

\textbf{(b)} As suggested by identities (\ref{sol35}) and (\ref{sol36}),
it is possible to simplify the statement of Theorem \ref{casfun}.
Since $T_{g_1}= T_{f_1}^{-1}$, we have\vspace*{-2pt}
%
\begin{equation} \label{xyenfuv}
(X,Y)= T_{g_1}(U,V) = \biggl( \log\biggl( \frac{1+U+V}{U+V}  \biggr), \log\biggl( \frac{1+U}{U}
\biggr) - \log\biggl( \frac{1+U+V}{U+V} \biggr) \biggr).\vspace*{-2pt}\vadjust{\goodbreak}
\end{equation}
As (\ref{xyenfuv}) shows, it is useful to introduce
%
\begin{equation} \label{uprimevprime}
(U^\prime, V^\prime) = \biggl( \frac{1+1/(U+V)}{ 1+1/U}, U+V\biggr).
\end{equation}
Obviously, the correspondence $(U,V)\mapsto(U^\prime, V^\prime)$ is
one-to-one:
%
\begin{equation} \label{uvenfuprimevprime}
(U,V) = \biggl( \frac{U^\prime V^\prime}{V^\prime+1 - U^\prime V^\prime},
\frac{V^\prime(V^\prime+1)(1-U^\prime)}{V^\prime+1 - U^\prime V^\prime}
\biggr).
\end{equation}
Furthermore, $(X,Y)$ can
be easily expressed in terms of $(U^\prime, V^\prime)$:
%
\begin{equation} \label{xyenfuprimevprime}
X=\log(1+1/V^\prime)\quad \mbox{and}\quad   Y= -\log U^\prime.
\end{equation}
Since it is easy to determine the density function of $\phi(\xi)$
knowing the density function of a~random variable $\xi$, where $\phi$
is differentiable and bijective, Theorem \ref{casfun} and its analog
related to $f=g_1$ (see Remark \ref{newrem}) are equivalent to Theorem
\ref{cas1} below.

\begin{theo} \label{cas1}
\textup{\textbf{(a)}} Let $U^\prime$ and $V^\prime$ be two positive and independent
random variables. The random variables $U$ and $V$ defined by
(\ref{uvenfuprimevprime}) are independent if only if there exist some
constants $a$, $b$, $c$ such that
%
\begin{equation} \label{uprimebetavprimek}
U^\prime\sim\operatorname{Beta} (a,b)\quad \mbox{and}\quad   V^\prime\sim
K^{(2)}(a+b,-b,c).
\end{equation}
If one of these equivalent
conditions holds, then
$U \sim K^{(2)}(a,b,p)$ and $V \sim\gamma(b,c)$.

\textup{\textbf{(b)}} Let $U$ and $V$ be two positive and independent random
variables. The random variables $U^\prime$ and $V^\prime$ defined by
(\ref{uprimevprime}) are independent if only if there exist some
constants~$a$, $b$, $c$ such that
%
\begin{equation} \label{ukvgamma}
U \sim K^{(2)}(a,b,c)\quad\mbox{and}\quad   V \sim\gamma(b,c).
\end{equation}
Under (\ref{ukvgamma}), $U^\prime\sim\operatorname{Beta} (a,b)$ and $V^\prime
\sim K^{(2)}(a+b,-b,c)$.
\end{theo}

We now formulate a simple consequence of Theorem \ref{cas1}.
\begin{theo} \label{nouveauthm}
For any $a,b,c>0$, the transformation $(u,v) \mapsto
(\frac{1+1/(u+v)}{ 1+1/u}, u+v)$ maps the probability
measure $K^{(2)}(a,b,c) \otimes\gamma(b,c)$ to the probability measure
$\operatorname{Beta} (a,b) \otimes K^{(2)}(a+b,-b,c).$ In particular,
%
\begin{equation} \label{cvkummer}
K^{(2)}(a,b,c) * \gamma(b,c)= K^{(2)}(a+b,-b,c).
\end{equation}
\end{theo}

\begin{rem} \label{rem14}
Note that (\ref{cvkummer}) may be regarded as an analog of
(\ref{convogig}).
\end{rem}

\subsection{\texorpdfstring{The case $f=f_\delta^*$}{The case f=f_delta^*}}
Recall that $f_\delta^{*}$ has been defined by (\ref{involution}).
Due to the form of $f_\delta^{*}$, a change of variables allows us to
simplify the search for independent random variables $X$ and $Y$ such
that the two components of $T_{f_\delta^{*}}(X,Y)$ are
independent.\vadjust{\goodbreak}

For any decreasing and
bijective function $f\dvtx (0,\infty) \rightarrow(0,\infty),$ we define
%
\begin{eqnarray} \label{fbarre}
\overline{f}(x)&=&\exp\{-f(-\log x)\},\qquad   x \in(0,1),\\ \label{nouvelletransf}
T_{f}^m(x,y)&=&\biggl( f(xy), \frac{f(x)}{f(xy)} \biggr),\qquad   x, y \in (0,1).
\end{eqnarray}

Observe that $\overline{f}$ is one-to-one and onto from $(0,1)$ to $(0,1)$,
$T_{f}^m$ is one-to-one and onto from $(0,1)^2$ to $(0,1)^2$ and
%
\begin{equation} \label{nouvelletransfinv}
( T_{f}^m )^{-1} = T_{f^{-1}}^m.
\end{equation}

\begin{defini} \label{defmult}
Let $X$ and $Y$ be two independent and $(0,1)$-valued random variables.
We say that a decreasing and bijective function $f\dvtx (0,1 ) \rightarrow
(0,1)$ is a multiplicative LWMY function with respect to $(X,Y)$ if the
random variables $U^m:=f(XY)$ and $V^m:=\frac{f(X)}{f(XY)}$ are
independent.
\end{defini}

\begin{rem} \label{remmult}
For any random vector $(X,Y)$ in $(0,\infty)^2$, we consider
$X^{\prime} = \mathrm{e}^{-X}$ and $Y^{\prime} = \mathrm{e}^{-Y}.$ Then, $f$ is an LWMY
function with respect to $(X,Y)$ if and only if $\overline{f}$ is a~multiplicative LWMY function with respect to $(X^{\prime},Y^{\prime})$.
\end{rem}

The change of variable $x^{\prime} = \mathrm{e}^{-x}$ is very
convenient since the function
%
\begin{equation} \label{fetoilebarre}
\phi_\delta(x):=\overline{f_\delta^{*}}(x)=\frac{1-x}{1+(\delta -1)x},\qquad
 x\in(0,1)
\end{equation}
is\vspace*{1pt} homographic.

Note that $\overline{f_\delta^{*}}\dvtx (0,1) \rightarrow(0,1)$ is
bijective, decreasing and equal to its inverse.
First, let us determine the distribution of the couple
$(X^{\prime},Y^{\prime})$ of random variables such that $\phi_\delta$
is a multiplicative LWMY function with respect to
$(X^{\prime},Y^{\prime})$.

For $ a, b, \alpha>0$ and $c\in\reel,$ consider the probability
measure
%
\begin{equation} \label{betaalpha}
\beta_\alpha(a,b;c)(\mathrm{d}x)=k_\alpha(a,b;c) x^{a-1} (1-x)^{b-1} (\alpha
x+1-x)^c \mathbf{1}_{(0,1)}(x)\,\mathrm{d}x.
\end{equation}
Note that if $c=0$, then
$\beta_\alpha(a,b;c)= \operatorname{Beta} (a,b)$.
\begin{theo} \label{theoremmult}
Let $X^\prime$ and $Y^\prime$ be two independent random variables
valued in $(0,1)$. Consider
\[
(U^m ,V^m) = T_{\phi_\delta}^m (X^\prime,Y^\prime) = \biggl( \frac{1-X^\prime
Y^\prime}{1+ (\delta-1) X^\prime Y^\prime},  \frac{1-X^\prime}{1+
(\delta-1) X^\prime} \frac{1+ (\delta-1) X^\prime Y^\prime}{1-X^\prime
Y^\prime} \biggr)
\]
for fixed $\delta>0$.

Then, $U^m$ and $V^m$ are independent if and only if there exist $a,
b, \lambda>0$ such that
%
\begin{equation} \label{loixprime}
X^\prime\sim\beta_\delta (a+b,\lambda;-\lambda-b),\qquad   Y^\prime \sim
\operatorname{Beta} (a,b).\
\end{equation}
If this condition holds, then
%
\begin{equation} \label{loium}
U^m \sim\beta_\delta (\lambda+b,a ;-a -b),\qquad   V^m \sim\operatorname{Beta}
(\lambda,b).
\end{equation}
\end{theo}

In the case $\delta=1$, Theorem \ref{theoremmult} takes a very simple
form.

\begin{propo}
Let $X^\prime$ and $Y^\prime$ be two independent random variables
valued in $(0,1)$. Then,
\[
U^m = 1-X^\prime Y^\prime,\qquad  V^m = \frac{1-X^\prime }{1 - X^\prime
Y^\prime}
\]
are independent if and only if
there exist $a, b, \lambda>0$ such that
\[
X^\prime\sim\operatorname{Beta} (a+b,\lambda)\quad  \mbox{and}\quad   Y^\prime\sim
\operatorname{Beta} (a,b).
\]
If one of these conditions holds, then $U^m \sim\operatorname{Beta} (\lambda
+b,a)$ and $V^m \sim\operatorname{Beta} (\lambda,b)$.
\end{propo}
\begin{rem}
When $ X^\prime\sim\operatorname{Beta} (a+b,\lambda)$ and $ Y^\prime
\sim\operatorname{Beta} (a,b),$ it can be proven that $U^m$ and $V^m$ are
independent using the well-known property that if $Z$ and $Z^\prime$
are independent with $Z \sim\gamma (a,1)$ and $Z^\prime \sim
\gamma (b,1),$ then $R:=\frac{Z}{Z+Z^\prime}$ and $Z+Z^\prime$ are
independent with $R \sim\operatorname{Beta} (a,b )$ and $Z+Z^\prime\sim
\gamma (a+b,1)$ (see, e.g., \cite{yorbeta}).
\end{rem}

According to Remark \ref{remmult}, $f_\delta^{*}$ is an LWMY function
with respect to $(X,Y)$ if and only if $\phi_\delta$ is a
multiplicative LWMY function with respect to $(X^{\prime},Y^{\prime}) =
(\mathrm{e}^{-X},\mathrm{e}^{-Y}) $. Therefore, a classical change of variables allows us
to deduce that Theorem \ref{theoremmult} is equivalent to Theorem
\ref{casfdelta} below.

\begin{theo} \label{casfdelta}
\textup{(1)} Consider two positive and independent random variables $X$ and $Y$.
The random variables $U=f_\delta^*(X+Y)$, $V=f_\delta^*(X)-f_\delta
^*(X+Y)$ are independent if and only if the densities of $Y$ and $X$
are, respectively,
%
\begin{eqnarray} \label{secondcase45}
p_Y(y)& = & \frac{\Gamma(a+b)}{\Gamma(a) \Gamma(b)} (1-\mathrm{e}^{-y})^{b-1}
\mathrm{e}^{-a y} \mathbf{1}_{ \{y>0\} }, \\[-2pt] \label{secondcase46}
p_X(x) & = &
k_\delta(a+b,\lambda, -\lambda-b) \mathrm{e}^{-(a+b ) x} (\delta \mathrm{e}^{-x}+1 -
\mathrm{e}^{-x})^{-\lambda-b }\nonumber\\ [-9pt]\\ [-9pt]
& & {}\times(1 - \mathrm{e}^{-x})^{\lambda-1 } \mathbf{1}_{x>0},\nonumber
\end{eqnarray}
where $a,b>0$, $\lambda\in\reel$ and $k_\delta(a+b,\lambda,
-\lambda-b)$ is the normalizing factor (see (\ref{betaalpha})). Thus,
$\mathrm{e}^{-Y}$ is $\operatorname{Beta}(a,b)$-distributed and $\mathrm{e}^{-X}$ is $\beta_\delta
(a+b, \lambda, -\lambda-b)$-distributed.

\textup{(2)} If \textup{(1)} holds, then the densities of $U$ and $V$ are, respectively,
%
\begin{eqnarray} \label{loiU}
p_U(u) & = & k_\delta(\lambda+b, a; -a -b) \mathrm{e}^{-u(\lambda+b)}
(1-\mathrm{e}^{-u})^{a-1} \nonumber\\ [-9pt]\\ [-9pt]
& &{} \times \bigl(1+(\delta-1)\mathrm{e}^{-u}\bigr)^{-a-b } \mathbf{1}_{u>0}, \nonumber\\[-2pt]\label{secondcase47}
 p_V(v) & = & \mathrm{e}^{-\lambda v} (1-\mathrm{e}^{-v})^{b-1}
\mathbf{1}_{v>0}.
\end{eqnarray}
\end{theo}

We omit the proof of Theorem \ref{casfdelta} since it is similar to
that of Theorem \ref{casfun}.\vadjust{\goodbreak}

\section{The set of all possible ``smooth'' LWMY functions}\label{sec3}
\setcounter{equation}{0}
The following theorem gives a functional
equation linking LWMY functions to the related densities.

\begin{theo} \label{theocns}
Let $X$ and $Y$ be two independent and positive random variables whose
densities $p_X$ and $p_Y$ are positive and twice differentiable.
Define $\phi_X= \log p_X$ and $\phi_Y=\log p_Y$.
Consider a decreasing function $f\dvtx (0,\infty) \mapsto(0,\infty)$,
three times differentiable. Then,
$f$ is a LWMY function with respect to $(X,Y)$
if and only if
%
\begin{eqnarray} \label{cns}
&&\phi_X^{\prime\prime}(x) - \phi_X^{\prime}(x) \frac{f^{\prime
\prime}(x)}{f^{\prime}(x)} + \phi_Y^{\prime\prime}(y) f^{\prime}(x) \biggl(
\frac {1}{f^{\prime}(x)} - \frac{1}{f^{\prime}(x+y)} \biggr)\nonumber\\
&&\quad{}+ \phi_Y^{\prime}(y) \frac{f^{\prime\prime}(x)}{f^{\prime}(x)} +
\frac{2(f^{\prime\prime}(x))^2 - f^{\prime\prime\prime}(x)
f^{\prime}(x)}{f^{\prime}(x)^2} =0,\qquad   x,y>0.
\end{eqnarray}
\end{theo}

\begin{pf}
Let $g=f^{-1}$ and $(U,V)=T_f(X,Y)$. By
formula (\ref{transfinverse}), $(X,Y)= T_g(U,V)$.
$X$ and~$Y$ being independent, the density of $(U,V)$ is
%
\begin{equation} \label{dens2}
p_{(U,V)} (u,v)=p_X\bigl(g(u+v)\bigr)p_Y\bigl(g(u)-g(u+v)\bigr) |J(u,v)| \mathbf{1}_{u,v>0},
\end{equation}
where $J$ is the Jacobian of the transformation $T_f$. We get
$|J(u,v)|= g^{\prime}(u+v)g^{\prime}(u)$,
and 
then
%
\begin{equation} \label{dens3}
p_{(U,V)} (u,v)=p_X\bigl(g(u+v)\bigr) p_Y\bigl(g(u)-g(u+v)\bigr)
g^{\prime}(u+v)g^{\prime}(u).
\end{equation}


The variables $U$ and $V$ are independent if and only if the function
$H=\log p_{(U,V)}$ satisfies $\frac{\partial^2 H}{\partial u \partial
v} =0.$ By equation (\ref{dens3}) we obtain
%
\begin{eqnarray} \label{deriv4}
\frac{\partial^2 H}{\partial u \,\partial v} & = & \phi_X^{\prime
\prime}(x) [g^{\prime}(f(x))]^2 +
\phi_X^{\prime}(x) g^{\prime\prime}(f(x)) \nonumber\\
& & {}-\phi_Y^{\prime\prime}(y) g^{\prime}(f(x)) \bigl[g^{\prime
}\bigl(f(x+y)\bigr) -g^{\prime}(f(x))\bigr] \nonumber\\
& &{} -\phi_Y^{\prime}(y) g^{\prime\prime}(f(x)) +\frac{ g^{\prime
\prime\prime} g^{\prime} - (g^{\prime\prime})^2}{(g^{\prime})^2}
(f(x)),
\end{eqnarray}
where $x= g(u+v)$ and $y=g(u)-g(u+v)$.
Differentiating three times the relation $g(f(x))=x$, we obtain
$ g^{\prime\prime}(f(x)) =
-\frac{f^{\prime\prime}(x)}{f^{\prime}(x)^3}$ and
$g^{\prime\prime\prime}(f(x)) =
-\frac{f^{\prime\prime\prime}(x)
f^{\prime}(x) - 3 f^{\prime\prime}(x)^2}{f^{\prime}(x)^5}.$ As a~result,
%
\begin{equation} \label{deriv7}
\frac{ g^{\prime\prime\prime} g^{\prime} - (g^{\prime
\prime})^2}{(g^{\prime})^2}
(f(x))
= \frac{ 2 f^{\prime\prime}(x)^2 - f^{\prime\prime\prime}(x)
f^{\prime}(x) }{f^{\prime}(x)^4}.
\end{equation}
%
Therefore, $\frac{\partial^2 H}{\partial u\, \partial v} =0$ leads to
(\ref{cns}).\
\end{pf}

We restrict ourselves to \textit{smooth} LWMY functions $f$, that is,
those satisfying
%
\begin{equation} \label{cond1}
f\dvtx(0,\infty)\rightarrow(0,\infty)  \mathrm{\ is\ bijective\ and \
decreasing,}
\end{equation}
%
\begin{equation} \label{cond2}
f  \mathrm{\ is\ three\ times\ differentiable,}
\end{equation}
%
\begin{equation} \label{serie}
F(x)= \sum_{n\geq1} a_nx^n\qquad  \forall x>0,
\end{equation}
where $F:=1/f^\prime$.

According to (\ref{cond1}), $f^{\prime}(0_+)=-\infty$. This implies
that $F(0_+)=0$ and explains why the series in (\ref{serie}) starts
with $n=1$.

The goal of this section is to prove half of Theorem \ref{theoclass}:
if $f$ is a smooth LWMY function, then $f$ belongs to one of the four
classes ${\mathcal F}_1,\ldots,{\mathcal F}_4$ introduced in Remark
\ref{remark1}. First, in Theorem~\ref{classificationF}, we
characterize all possible functions $F$. Second, we determine the
associated functions $f$ (see Theorem \ref{classificationFbis}).

\begin{theo} \label{classificationF}
Suppose that $f$ is a smooth LWMY function and the assumptions of
Theorem \ref{theocns} are satisfied.
\begin{enumerate}
\item If $F^{\prime}(0_+)=0$, then $a_2 <0$ and
%
\begin{equation} \label{premiercas}
F(x)= \cases{
\displaystyle\frac{a_2^2}{6a_4} \biggl( \cosh\biggl( x \sqrt{ \frac{12a_4}{a_2}
} \biggr)-1 \biggr), &\quad if $a_4 < 0$, \cr
a_2x^2, &\quad  otherwise.
}
\end{equation}
\item If $F^{\prime}(0_+)\neq0$, then

\begin{equation} \label{secondcas}
F(x)= \cases{
\displaystyle\frac{a_1a_2}{3a_3} \biggl[ \cosh\biggl( x \sqrt{\frac{6a_3}{a_1}} \biggr)-1 \biggr]\cr
\quad\displaystyle{} + a_1
\sqrt{\frac{a_1}{6a_3}} \sinh\biggl( x \sqrt{\frac
{6a_3}{a_1}} \biggr), &\quad if $a_1 a_3 > 0$, \cr
a_1x + a_2x^2, &\quad otherwise.
}
\end{equation}
\end{enumerate}
\end{theo}

\begin{rem} \label{casyor}
Unsurprisingly, the case $F(x)=a_2x^2$ corresponds to
$f(x)=-\frac{1}{a_2} \frac{1}{x}$, that is, the case considered by
Matsumoto and Yor, and Letac and Wesolowski. Thus, under stronger
assumptions, we retrieve the result of Letac and Wesolowski. Indeed,
writing the functional equation of Theorem \ref{cns} with $f\dvtx
x\mapsto1/x$ gives
\[
\phi_X^{\prime\prime}(x) + \frac{2}{x}\phi_X^{\prime}(x)
+ \phi_Y^{\prime\prime}(y) \frac{1}{x^2} \bigl(x^2-(x+y)^2\bigr)
- \frac{2}{x} \phi_Y^{\prime}(y) +\frac{2}{x^2}=0.
\]
We then solve this differential equation and find that the laws of $X$
and $Y$ are necessarily $\operatorname{GIG}$ and gamma, respectively.
We omit the details.
\end{rem}

Throughout this subsection, we suppose that $f$ satisfies
(\ref{cond1})--(\ref{serie}) and that the assumptions of Theorem
\ref{theocns} are fulfilled. To simplify the statement of results
below, we do not repeat these conditions.

Recall that $\phi_Y$ is the logarithm of the density of $Y$. Let us
introduce
%
\begin{equation} \label{loi13}
h:=\phi_Y^{\prime}.\vadjust{\goodbreak}
\end{equation}

\begin{lemma} \label{premierlemme}
%
\begin{enumerate}
\item There exists a function $\lambda\dvtx  (0,\infty) \rightarrow \reel$
such that
%
\begin{equation} \label{grantef15}
F(x+y)=\frac{\lambda(x)-h(y) F^{\prime}(x)}{h^{\prime}(y)}+ F(x) .
\end{equation}
\item $F$ satisfies
%
\begin{equation} \label{grantef16}
F(y)=\frac{\lambda(0_+)-h(y)F^{\prime}(0_+)}{h^{\prime}(y)}.
\end{equation}
\end{enumerate}
\end{lemma}

\begin{rem}
Suppose that we have been able to determine $F$. Then,
$h=\phi_Y^{\prime}$ solves the linear ordinary differential equation
(\ref{grantef16}) and can therefore be determined. The remaining
function $\phi_X$ is obtained by solving equation (\ref{cns}).
\end{rem}

\begin{pf*}{Proof of Lemma \ref{premierlemme}}
Using (\ref{loi13}) and $F=1/f^{\prime}$ in equation (\ref{cns}), we
obtain
%
\[
c(x)=h(y)\frac{F^{\prime}(x)}{F(x)}+h^{\prime}(y)
\frac{1}{F(x)}\bigl(F(x+y) -F(x)\bigr),
\]
where $c(x)$ depends only on $x$.
Multiplying both sides by $F(x)$ and taking the $y$-derivative leads
to
\[
0= F^{\prime}(x)
h^{\prime}(y)+ [F(x+y)-F(x)]h^{\prime\prime}(y)
+ h^{\prime}(y)
F^{\prime}(x+y).
\]

Fix $x>0$. Then, $\theta(y):= F(x+y)$ is a solution of the
differential equation in $y$

\begin{equation} \label{grantef17}
0= F^{\prime}(x) h^{\prime}(y)+ \bigl(\theta(y)-F(x)\bigr)h^{\prime\prime}(y)
+
h^{\prime}(y) \theta^{\prime}(y).
\end{equation}
A solution of the related homogeneous equation in $y$ is
$\frac{\rho}{h^{\prime}(y)}$, where $\rho$ is a constant. It is
easy to
prove that $y\mapsto - F^{\prime}(x) h(y) +F(x)h^{\prime}(y)$ solves
(\ref{grantef17}). Thus, the general solution of (\ref{grantef17}) is
\[
\theta(y) = \frac{1}{h^{\prime}(y)} [ \lambda(x) - F^{\prime }(x) h(y)
+F(x)h^{\prime}(y) ] .
\]
Since $\theta(y)= F(x+y)$, (\ref{grantef15}) follows.

According to (\ref{serie}), $F(0_+)$ and $F^{\prime}(0_+)$ exist.
Therefore, taking the limit $x \to0_+$ in~(\ref{grantef15}) implies
both the existence of $\lambda(0_+)$ and relation (\ref{grantef16}).
\end{pf*}

The following lemma shows that the function $F$ (and thus $f$) solves
a self-contained equation in which $h$, and thereby the densities of
$X$ and $Y$, are not involved.

\begin{lemma}
$F$ solves the delay equations
%
\begin{eqnarray}\label{grantef18}
F(x+y) & = &\frac{F(y)[ \lambda(x) -h(y)F^{\prime}(x) ]}{ \lambda
(0_+) -h(y)F^{\prime}(0_+) } +F(x)\qquad  (x,y > 0), \\\label{grantef19}
F^{\prime}(x+y) & = & \frac{F^{\prime}(y)+
F^{\prime}(0_+) }{F(y)} [ F(x+y)-F(x)] -F^{\prime}(x)\qquad(x,y > 0).
\end{eqnarray}
\end{lemma}

\begin{pf}
By (\ref{grantef16}), we have
\[
h^{\prime}(y)= \frac{\lambda(0_+)-h(y)F^{\prime}(0_+)}{F(y)}.
\]
Equation (\ref{grantef18}) then follows by rewriting equation (\ref
{grantef15}) and replacing $h^{\prime}(y)$ with the expression above.

We differentiate (\ref{grantef18}) in $y$ and use the fact that $
\lambda(0_+)-h(y)F^{\prime}(0_+)= h^{\prime}(y)F(y)$
to obtain~

\[
F^{\prime}(x+y) = [F^ \prime(y)+ F^ \prime(0_+)]
\frac{\lambda(x)
-h(y)F^{\prime}(x)}{F(y) h^{\prime}(y)}- F^{\prime}(x).
\]
By (\ref{grantef15}), we have
$\frac{\lambda(x) -h(y)F^{\prime}(x)}{F(y) h^{\prime}(y)} =
\frac{F(x+y)-F(x)}{F(y)}$ and this gives (\ref{grantef19}).
\end{pf}

\begin{rem} \label{delay}
We can see (\ref{grantef19}) as a scalar neutral delay
differential equation. Indeed, set $t=x+y$ and consider $y>0$ as a
fixed parameter. Then, (\ref{grantef19}) becomes
%
\begin{equation} \label{eqdelay}
F^\prime(t)=a\bigl(F(t)-F(t-y)\bigr)-F^\prime(t-y),\qquad     t\geq y,
\end{equation}
where $a:=\frac{F^\prime(y)+F^\prime(0_+)}{F(y)}$. Replacing $F(t)$
in (\ref{eqdelay}) with $\mathrm{e}^{at}G(t)$ leads to
%
\begin{equation} \label{eqdelaybis}
G^\prime(t)+ \mathrm{e}^{-ay}G^\prime(t-y)+2a \mathrm{e}^{-ay}G(t-y)=0,\qquad    t\geq y.
\end{equation}
Equation (\ref{eqdelaybis}) is called a \textit{neutral delay
differential equation} (see, e.g., Section 6.1, in~\cite{gyor}). These equations have been intensively studied, but the authors
have only focused on the asymptotic behavior of the solution as $t \to
\infty$. Unfortunately, these results do not help to solve explicitly
either (\ref{grantef19}) or (\ref{eqdelaybis}).
\end{rem}
%

\begin{lemma} \label{combinatories}
For all integers $k\geq0$ and $l\geq1$, we have
%
\begin{eqnarray} \label{comb20}
\sum_{m=0}^{l-1}(l-2m+1)C_{l-m+1+k}^k a_{l-m+1+k} a_m & =&
(l-2)(k+1)a_{k+1}a_l+a_1a_{l+k} C_{l+k}^k ,\\ \label{comb21}
C_{k+3}^k a_{k+3}a_1 & =& (k+1)a_{k+1}a_3 ,
\end{eqnarray}
\vspace*{-17pt}
\begin{equation} \label{comb22}
2C_{k+4}^k a_{k+4}a_1 + C_{k+3}^k a_{k+3}a_2 - C_{k+2}^k a_{k+2}a_3
- 2 (k+1) a_{k+1}a_4 =0 ,
\end{equation}
where $C_n^p=\frac{n!}{(n-p)!p!}$.
\end{lemma}

\begin{pf} Obviously, the equation (\ref{grantef19}) is equivalent to
%
\begin{eqnarray} \label{comb23}
F^{\prime}(x+y) F(y) &=& F^{\prime}(y) F(x+y)- F^{\prime}(y) F(x)\nonumber\\ [-8pt]\\ [-8pt]
&&{} -F(y)
F^{\prime}(x) + F^{\prime}(0_+) F(x+y) - F^{\prime}(0_+) F(x).\nonumber
\end{eqnarray}
Using the asymptotic expansion (\ref{serie}) of $F,$ we can develop
each term in (\ref{comb23}) as a~series with respect to $x$ and $y$.
Then, identifying the series on the right-hand side and the left-hand
side, we get (\ref{comb20})--(\ref{comb22}).
The details are provided in the \hyperref[appendix]{Appendix}.
\end{pf}

\begin{pf*}{Proof of Theorem \ref{classificationF}}
We will only prove item 1;
the proof of item 2 is similar.

Since $a_1=F^{\prime}(0_+)=0$, we necessarily have $a_2\neq0$. Indeed,
if $a_2=0,$ then, by (\ref{comb22}) with $k=1$, we would have
$-3a_3^2-4a_2a_4=0$, that is, $a_3=0$. Again using (\ref{comb22})
with $k=3$ would imply that $a_4=0$ and finally that $a_k=0$ for every
$k\geq0$, which is a contradiction because, by definition,
$F=1/f^{\prime}$ does not vanish.

So, we have $a_1=0$ and $a_2\neq0$. Equation (\ref{comb21}) with $k=1$
reads $4a_4a_1 =2a_2a_3,$ which implies that $a_3=0$. Applying
(\ref{comb21}) to $k=2n$ provides, by induction on $n$, $a_{2n+1}=0$
for every $n\geq0$.

Therefore, equation (\ref{comb22}) reduces to
$(k+3)(k+2)(k+1)a_{k+3}a_2 = 12 (k+1) a_{k+1}a_4,   k\geq0$, that is,
$a_{k+3}= \frac{12a_4}{a_2} \frac{1}{(k+3)(k+2)}a_{k+1}.$ This leads to
%
\begin{equation} \label{comb33}
a_{2k} =\biggl( \frac{12a_4}{a_2} \biggr) ^{k-1} \frac{2}{(2k)!} a_2,\qquad     k\geq1.
\end{equation}
Then, $F(x)=a_2x^2$ if $a_4=0,$ and if
$a_4\neq0,$ we have
\[
F(x) = \sum_{k\geq1} \biggl( \frac{12a_4}{a_2} \biggr) ^{k-1} \frac {2}{(2k)!} a_2
x^{2k}.
\]
If $a_4 a_2<0$, then $F(x)=\frac{a_2^2}{6a_4} [ \cos( x \sqrt{ \frac
{-12a_4}{a_2} } ) -1 ].$ This implies $F(2\uppi\sqrt{ \frac{-12a_4}{a_2}
})=0,$ which is impossible since $F(x)=1/f^\prime(x)<0$. Consequently,
\[
F(x)=\frac{a_2^2}{6a_4} \biggl[ \cosh\biggl( x \sqrt{ \frac {12a_4}{a_2} } \biggr) -1 \biggr].
\]
\upqed\end{pf*}

Now, in each case of Theorem \ref{classificationF}, we compute the
function $f$ associated with $F$ via the relation $F=1/f^{\prime}$. We
do not detail the calculations since they reduce to getting a good
primitive of $1/F$. Recall that we restrict ourselves to functions $f$
satisfying (\ref{cond1})--(\ref{serie}) and work under the assumptions
of Theorem \ref{cns}.\

\begin{theo} \label{classificationFbis}
\begin{enumerate}
\item If $F(x)=a_2 x^2$, then $f(x)=\frac{1}{a_2 x}$.
\item If $F(x)=
\alpha( \cosh\beta x -1), \alpha, \beta>0 $, then
$f(x)=\frac{2}{\alpha\beta}f_1(\beta x)$.
\item If $F(x)= a_1x+a_2x^2$, then $f(x)=
-\frac{1}{a_1}g_1(\frac{a_2}{a_1}x)$.
\item If
\[ F(x)= \frac{a_1a_2}{3a_3} \biggl[ \cosh\biggl( x
\sqrt{\frac{6a_3}{a_1}} \biggr)-1 \biggr] + a_1 \sqrt{\frac{a_1}{6a_3}} \sinh\biggl( x
\sqrt{\frac{6a_3}{a_1}} \biggr),
\]
then
\[
f(x)=-\frac{1}{\beta\gamma} \log\biggl( \frac{e^{\beta x}+\delta -1}{e^{\beta
x} -1} \biggr),
\]
where $\alpha= \frac{a_1a_2}{3a_3}$, $\beta= \sqrt{\frac
{6a_3}{a_1}}$ and
$\gamma= a_1 \sqrt{\frac{a_1}{6a_3}}$.
\end{enumerate}
\end{theo}

\section{\texorpdfstring{Proof of Theorem \protect\ref{casfun}}{Proof of Theorem 2.4}}\label{sec4}

Recall that $\phi_Y= \log p_Y$, $h=\phi_Y^{\prime}$ and
$F^{\prime}(0_+)=0$. It is easy to deduce from (\ref{grantef16}) that
there exist constants $\lambda$ and $c_1$ such that $h(y) = \lambda
f(y)+c_1$, that is, $h(y) = \frac{\lambda \mathrm{e}^y}{\mathrm{e}^y-1} +c_1 -\lambda$.
This implies the existence of a constant $d$ such that $\phi_Y (y)=
\lambda\log(\mathrm{e}^y-1)+ (c_1 -\lambda)y +d$. Setting $M=\mathrm{e}^d$, we have, by
integration, for all $y>0$,
%
\begin{equation} \label{sol38}
p_Y(y) = M(1-\mathrm{e}^{-y})^\lambda \mathrm{e}^{c_1 y}.
\end{equation}
To give more information on the normalizing constant $M$, we observe,
for $a=-c_1$ and $b=\lambda+1$, that
\[
\int_0^{\infty} M(1-\mathrm{e}^{-y})^{b-1} \mathrm{e}^{-a y}\,\mathrm{d}y = M \int_0^1 (1-u)^{b-1}
u^{a-1}\,\mathrm{d}u,
\]
which implies that $a>0$, $b>0$ and $M=\frac{\Gamma(a+b)}{\Gamma(a)
\Gamma(b)}$. This\vspace*{1pt} proves (\ref{sol35}).

To find the density of $X,$ we return to equation (\ref{cns}) and
compute each of its terms.

We have
$f^{\prime} (x) = \frac{-\mathrm{e}^x}{(\mathrm{e}^x-1)^2}$,
$f^{\prime\prime} (x) = \frac{\mathrm{e}^{2x} + \mathrm{e}^x}{(\mathrm{e}^x-1)^3}$ and
$f^{\prime
\prime\prime} (x)= -\frac{\mathrm{e}^{3x}+4\mathrm{e}^{2x} + \mathrm{e}^x}{(\mathrm{e}^x-1)^4} $ so that
$ \frac{ f^{\prime} (x)}{f^{\prime} (x+y)} = \frac{\mathrm{e}^{-y}
(\mathrm{e}^{x+y}-1)^2}{(\mathrm{e}^x-1)^2} $ and $ \frac{ f^{\prime\prime}
(x)}{f^{\prime} (x)} = -\frac{\mathrm{e}^x+1}{\mathrm{e}^x-1}. $ Calculations yield
%
\begin{equation} \label{sol41}
\frac{2(f^{\prime\prime}(x))^2 - f^{\prime\prime\prime}(x)
f^{\prime}(x)}{f^{\prime}(x)^2} =\frac{\mathrm{e}^{2x}+1}{(\mathrm{e}^x-1)^2}.
\end{equation}
Moreover,
%
\begin{equation} \label{sol42}
-\phi_Y^{\prime}(y) \frac{ f^{\prime\prime} (x)}{f^{\prime} (x)} +
\phi_Y^{\prime\prime}(y) \biggl( \frac{f^{\prime}(x)}{f^{\prime}(x+y)}-1 \biggr) =
\frac {(c_1-\lambda )\mathrm{e}^{2x} -c_1}{(\mathrm{e}^x-1)^2}.
\end{equation}
Equation (\ref{cns}) can then be written, using (\ref{sol41}) and
(\ref{sol42}),
\[
\phi_X^{\prime\prime}(x) + \frac{\mathrm{e}^x+1}{\mathrm{e}^x-1} \phi_X^{\prime}(x)
= \frac{(c_1-\lambda-1)\mathrm{e}^{2x} -c_1-1}{(\mathrm{e}^x-1)^2}.
\]
Then,
$h_0:=\phi_X^{\prime}$ solves
%
\begin{equation} \label{sol43}
h_0^{\prime}(x) + \frac{\mathrm{e}^x+1}{\mathrm{e}^x-1} h_0(x)= \frac{(c_1-\lambda
-1)\mathrm{e}^{2x} -c_1-1}{(\mathrm{e}^x-1)^2}.
\end{equation}
Note that $x \mapsto\frac{K}{4 \sinh^2(x/2)}$ solves (\ref{sol43})
with the right-hand side equal to 0, and $x \mapsto\frac{
(c_1-\lambda
-1)\mathrm{e}^x +(c_1+1)\mathrm{e}^{-x}}{4 \sinh^2(x/2)}$ is a particular solution of
(\ref{sol43}). Therefore, the solution of (\ref{sol43}) is
\[
h(x)= \frac{ (c_1-\lambda-1)\mathrm{e}^x +(c_1+1)\mathrm{e}^{-x}+ K }{4\sinh^2(x/2)}
\]
for some constant $K$. This implies that
\[
\phi_X^{\prime}(x)
= c_1+1 +\frac{(2c_1-\lambda+K)\mathrm{e}^x}{(\mathrm{e}^x-1)^2}-\frac{(\lambda
+2)\mathrm{e}^x}{\mathrm{e}^x-1}.
\]
As a consequence, there exists a constant $\delta$
such that
\[
\phi_X(x)= (c_1+1)x -\frac{(2c_1-\lambda+K)\mathrm{e}^x}{\mathrm{e}^x-1}-(\lambda
+2)\log(\mathrm{e}^x-1)+\delta.
\]
Thus, $p_X(x)= N \mathrm{e}^{(c_1+1 ) x} (\mathrm{e}^x-1)^{-\lambda-2} \exp(
-\frac{2c_1-\lambda+ K}{\mathrm{e}^x-1} ) \mathbf{1}_{ \{x>0\} }.$ Recall that
$a=-c_1$ and $b=\lambda+1$. With $c= 2c_1-\lambda+ K$, we get
(\ref{sol36}). More information on the constant $N$ is obtained by
observing that if we set $V^\prime= f_1(X)=\frac{1}{\mathrm{e}^X-1}$, then the
density of $V^\prime$ is
\[
f_{V^\prime}(w)=N (w+1)^{-a}w^{a+b-1} \exp\{ -cw \} \mathbf{1}_{ \{w>0\} },
\]
that is, the law of $V^\prime$ is $K^{(2)} (a+b, -b,c)$ (see equation
(\ref{kummer})).

We have $g_1^{\prime}(u)= -\frac{1}{u(u+1)}. $ A computation of a
Jacobian,
together with (\ref{sol35}) and~(\ref{sol36}), implies, for $u,v>0$,
that
\begin{eqnarray*} p_{(U,V)} (u,v) & = & p_X \biggl( \log
\biggl[\frac{u+v+1}{u+v}\biggr] \biggr) p_Y \biggl( \log\biggl[ \frac{(u+1)(u+v)}{u(u+v+1)} \biggr] \biggr)\\
& &{} \times\frac{1}{u(u+1)(u+v)(u+v+1)}.
\end{eqnarray*}
We then get that $p_{(U,V)} (u,v)$ is the product of a function of $u$
and a function of $v,$ and this gives item 2 of Theorem \ref{casfun}.


\begin{appendix}\label{appendix}
\section*{Appendix}
\setcounter{section}{5}
\renewcommand{\theequation}{\arabic{section}.\arabic{equation}}
\setcounter{equation}{4}
\begin{pf*}{Proof of Lemma \ref{combinatories}}
We have
\[
F^{\prime}(x+y) F(y) =
\sum_{k\geq0} x^k \sum_{m\geq0, n\geq1+k} n a_n a_m C_{n-1}^k
y^{n+m-1-k}.
\]
Setting $l=m+n-1-k$ for fixed $m$ gives
%
\begin{equation} \label{comb24}
F^{\prime}(x+y) F(y) = \sum_{k\geq0, l\geq0} x^k y^{l} \sum_{m=
0}^l (l-m+1+k) C_{l-m+k}^k a_{l-m+1+k} a_m.
\end{equation}
By the same method, we have
%
\begin{equation} \label{comb25}
F^{\prime}(y) F(x+y) = \sum_{k\geq0, l\geq0} x^k y^{l} \Biggl( \sum_{m=
0}^{l+1} m C_{l-m+k+1}^k a_{l-m+1+k} a_m \Biggr).
\end{equation}
As for the two other terms of (\ref{comb23}), we get
%
\begin{eqnarray} \label{comb26}
F^{\prime}(y) F(x) & = & \sum_{k\geq0, l\geq0} a_k a_{l+1} (l+1)
x^k y^l, \\\label{comb27}
 F^{\prime}(x) F(y) & = & \sum_{k\geq0, l\geq0}
a_{k+1} a_{l} (k+1) x^k y^l.
\end{eqnarray}
Consequently,
%
\begin{eqnarray} \label{comb28}
F^{\prime}(0_+) F(x+y) &=& a_1 \sum_{n \geq0} a_n (x+y)^{n}
= a_1 \sum_{k,l\geq0} a_{l+k} C_{l+k}^k x^k y^{l},\\ \label{comb29}
F^{\prime}(0_+) F(x) &=& a_1 \sum_{k \geq0} a_k x^{k}.
\end{eqnarray}

Identifying the coefficient of $x^ky^l$ in (\ref{comb23}) and using
(\ref{comb24})--(\ref{comb29}), we have, for $k\geq0$ and $ l \geq0$,
%
\begin{eqnarray}\label{comb30}
\sum_{m=0}^{l}(l-m+1+k)C_{l-m+k}^k a_{l-m+1+k} a_m & = & -(l+1)a_k
a_{l+1}-(k+1)a_{k+1}a_l\nonumber\\
& &{} +\sum_{m= 0}^{l+1} m C_{l-m+k+1}^k a_{l-m+1+k} a_m \\
& &{} + a_1a_{l+k} C_{l+k}^k - a_1 a_k 1_{l=0}.\nonumber
\end{eqnarray}
Note that if $l=0$, then both sides of (\ref{comb30}) vanish.
Therefore, we may suppose in the sequel that $l\geq1$.

For $m=l+1,$ we have $m C_{l-m+k+1}^k a_{l-m+1+k} a_m = (l+1)a_k
a_{l+1}$. Thus, equation (\ref{comb30}) reads
%
\begin{eqnarray}\label{comb31}
&&\sum_{m=0}^{l}(l-m+1+k)C_{l-m+k}^k a_{l-m+1+k} a_m\nonumber \\
&&\quad =
-(k+1)a_{k+1}a_l+ \sum_{m= 0}^{l} m C_{l-m+k+1}^k a_{l-m+1+k}
a_m\\
& &\qquad{} +a_1a_{l+k} C_{l+k}^k.\nonumber
\end{eqnarray}
However, via a calculation involving the definition, we find that
\[
(l-m+1+k)C_{l-m+k}^k - m C_{l-m+1+k}^k= (l-2m+1) C_{l-m+1+k}^k,
\]
so equation (\ref{comb31}) is equivalent to
%
\begin{equation} \label{comb32}
\sum_{m=0}^{l}(l-2m+1)C_{l-m+1+k}^k a_{l-m+1+k} a_m = -(k+1)a_{k+1}a_l+
a_1a_{l+k} C_{l+k}^k.
\end{equation}

For $m=l$, we have $(l-2m+1)C_{l-m+1+k}^k a_{l-m+1+k} a_m =(1-l) (k+1)
a_{k+1} a_l$. Consequently, equation (\ref{comb32}) may be written
as
\[
\sum_{m=0}^{l-1}(l-2m+1)C_{l-m+1+k}^k a_{l-m+1+k} a_m - (l-1) (k+1)
a_{k+1} a_l = -(k+1)a_{k+1}a_l+ a_1a_{l+k} C_{l+k}^l,
\]
which implies (\ref{comb20}).

(\ref{comb21}) and (\ref{comb22}) follow by applying (\ref{comb20}) to
$l=3$ and $l=4,$ respectively.
\end{pf*}
\end{appendix}

\section*{Acknowledgements} We are grateful to G.~Letac for helpful
discussions about this work, and to a referee whose comments led to an
improvement of the paper.

\printhistory

\end{document}